\documentclass[10pt]{amsart}

\usepackage{hyperref}

\usepackage{amssymb}
\usepackage{amsmath}
\usepackage{mathabx}

\let\newpf\proof \let\proof\relax

\def\bm{\begin{matrix}}
\def\em{\end{matrix}}

\newcommand{\bt}{\begin{thm}}
\newcommand{\et}{\end{thm}}

\newcommand{\bl}{\begin{lemma}}
\newcommand{\el}{\end{lemma}}

\newcommand{\beq}{\begin{eqnarray}}
\newcommand{\eeq}{\end{eqnarray}}

\def\be{\begin{equation}}
\def\ee{\end{equation}}

\def\ba{{\begin{align}}}
\def\ea{{\end{align}}}

\def\0{{\mathbf 0}}

\newtheorem{thm}{Theorem}[section]
\newtheorem{cor}[thm]{Corollary}

\newtheorem{lemma}[thm]{Lemma}

\theoremstyle{remark}

\numberwithin{equation}{section}

\def \bn {\hfill \\ \smallskip\noindent}

\theoremstyle{definition}

\newtheorem{definition}{Definition}[section]
\def\proof{\bn {\bf Proof.} }

\def\note#1
{\marginpar

{\tiny $\leftarrow$
\par
\hfuzz=20pt \hbadness=9000 \hyphenpenalty=-100 \exhyphenpenalty=-100
\pretolerance=-1 \tolerance=9999 \doublehyphendemerits=-100000
\finalhyphendemerits=-100000 \baselineskip=6pt
#1}\hfuzz=1pt}

\newcommand{\dist}{\operatorname{dist}}

\newcommand{\supp}{\operatorname{supp}}

\newcommand{\R}{{\mathbb R}}

\newcommand{\Z}{{\mathbb Z}}

\def\B0{{\bold{0}}}

\catcode`\@=12

\def\Empty{}
\newcommand\oplabel[1]{
  \def\OpArg{#1} \ifx \OpArg\Empty {} \else
  	\label{#1}
  \fi}

\newcommand{\comm}[1]{}
\newcommand{\comment}[1]{}

\begin{document}

\title{Shnol's theorem and the spectrum of long range operators}

\author{Rui Han}

\address{Department of Mathematics, University of California, Irvine}

\email{rhan2@uci.edu}

\thanks{}

\begin{abstract}
We extend some basic results known for finite range operators to long range operators with off-diagonal decay. Namely, we prove an analogue of Shnol's theorem. We also establish the connection between the almost sure spectrum of long range random operators and the spectra of deterministic periodic operators.
\end{abstract}

\maketitle

\section{Introduction}
Long range operators on $l^2(\Z^d)$ arise naturally from the discrete Laplacian on half-space $\Z^{d+1}_{+}$ after dimension reduction (see e.g. 
\cite{G, JM2,JMP}). 
However, compared to Schr\"odinger operators or finite range operators, little has been studied.
We refer the readers to \cite{ASV,SSV} for studies of some special models of long range operators and the integrated density of states.
Our goal is to extend some basic results known for finite range operators to long range setting with off-diagonal decay. 
The first part of this note concerns a generalization of Shnol's theorem to the case of a long range normal operator and a long range self-adjoint operator with unbounded potentials.
The second part provides a connection between the almost sure spectrum of long range self-adjoint random operators and the spectra of deterministic periodic operators.
We hope these could serve as ready-to-use tools in the future for people who study discrete operators. 

The classical Shnol's theorem for Schr\"odinger operators is well-known, see \cite{schnol, berez, Simon, SimonBAMS, book} for the continuum case, from which the discrete version could be derived. 
It asserts that any spectral measure gives full weight to the set of energies with polynomially bounded generalized eigenfunctions, moreover, the spectrum is the closure of this set. 
It has many applications, for example, relating the spectrum to non-uniform hyperbolicity of the corresponding cocycle \cite{Johnson,HP,Marx,DFLY16}, as well as providing a priori estimate which turns out to be crucial in the proofs of the almost-localization and localization. 
So far there have been various generalizations of the Shnol's theorem (see e.g. \cite{Abel, K, BLS, FLW, DFLY16,LT}). 
In this note, we generalize this result to the long range case.

Let us consider long range normal operators with polynomially decaying off-diagonal terms.

\begin{align}\label{defnormalH}
(Hu)(n)=\sum_{j\in \Z^d}a_{j}^{n}u(n-j),
\end{align}
with
\begin{align}\label{adecay}
|a_j^m|\leq C(1+\|j\|)^{-r}\ \mathrm{for}\ \mathrm{any}\ j, m\in \Z^d,
\end{align}
where $\|j\|=(\sum_{k=1}^d j_k^2)^{\frac{1}{2}}$ is the Euclidean norm of $j=(j_1,j_2,...,j_d)$, $C>0$ is a constant and $r>d/2$ is a constant whose value will be specified later.
Note that $r>d/2$ ensures that $\{a_{j}^m\}_{j\in \mathbb{Z}^d}\in l^2(\mathbb{Z}^d)$ with uniformly bounded $l^2$-norm in $m$.

We assume
\begin{align}\label{eqnormal}
\sum_{j\in \Z^d} a_j^n\overline{a_{j-m}^{n-m}}=\sum_{j\in \Z^d}\overline{a_{-j}^{n-j}}a^{n-j}_{m-j}\ \ \mathrm{for}\ \mathrm{any}\ m,n\in \Z^d, 
\end{align}
to ensure $H$ is a normal operator.

We introduce generalized eigenfunctions.
\begin{definition} ({\it $\epsilon$-generalized eigenvalue/eigenfunction})
Let $\epsilon$ be a positive number. An energy $z$ of $H$ is called an $\epsilon$-generalized eigenvalue if there is a formal solution to the equation $H\phi_z=z\phi_z$, with $\phi_z(0)=1$ and $|\phi_z(n)|\leq C^\prime (1+\|n\|)^{\frac{d}{2}+\epsilon}$ for some constant $C^\prime$.
\end{definition}

We denote the set of $\epsilon$-generalized eigenvalues by $G_\epsilon$. 
The spectral measures of $H$ are defined by $\mu_{n,m}(B)=(e_n, \chi_{B}(H)e_m)$, for any Borel sets $B\subset \mathbb{C}$ and $m,n\in \Z^d$. 
Here $\chi_B$ denotes the characteristic function of $B$. 
We denote $\mu_{n,n}:=\mu_n$. 
Let $\mu=\sum_{n\in \Z^d}\lambda_n\mu_n$, where $\lambda_n= (1+\|n\|)^{-d-2\epsilon}/(\sum_{n\in \Z^d} (1+\|n\|)^{-d-2\epsilon})$. Clearly, any spectral measure is absolutely continuous with respect to $\mu$.

Our main theorem is:
\begin{thm}\label{main}Let $H$ be a normal operator defined as in (\ref{defnormalH}) with $r>2d$. Then for any $0<\epsilon<r-2d$, we have the following:\\
$(a).\ \ G_{\epsilon}\subseteq\sigma(H)$,\\
$(b).\ \ \mu(\sigma(H)\setminus G_{\epsilon})=0$,\\
$(c).\ \ \overline{G_{\epsilon}}=\sigma(H)$.
\end{thm}
The proof of part (a) relies on Lemma \ref{key}, which is inspired by those for the Schr\"odinger case \cite{schnol,SimonBAMS}. 
We will discuss both in Section 2. The proofs of (b), (c) are standard. We will include them in the appendix for completeness.

\

Switching our attention to the long range self-adjoint case, we are able to cover operators with unbounded potentials $a_0^n$.
Let $H$ be a self-adjoint operator defined as follows.
\begin{align}\label{defsaH}
(Hu)(n)=\sum_{j\in \Z^d\setminus \{0\}}a_{j}^{n}u(n-j)+a_0^n u(n),
\end{align}
with
\begin{align}\label{adecay}
|a_j^m|\leq C(1+\|j\|)^{-r}\ \mathrm{for}\ \mathrm{any}\ j\neq 0, m\in \Z^d,\ \ \mathrm{and}\ \ a_{j}^n=\overline{a_{-j}^{n-j}}.
\end{align}

Note that the polynomial decay condition (\ref{adecay}) does not involve $j=0$, allowing us to take unbounded potentials into account.
Following exactly the same line of the proof of Theorem \ref{main}, we have Shnol's theorem for long range self-adjoint operators (with unbounded potentials).
\begin{thm}\label{samain}
Let $H$ be a self-adjoint operator defined as in (\ref{defsaH}) with $r>2d$. 
Then for any $0<\epsilon<r-2d$, the conclusions of Theorem \ref{main} hold.
\end{thm}

The localization property of long range self-adjoint operators with polynomially decaying off-diagonal terms has been studied in the random potential case, when the potentials $a_0^n=V_{\omega}(n)$ are i.i.d. random variables. 
Some known results are pure point spectrum in the large disorder and high energy regime when $r>d$ \cite{AM}, purely singular spectrum for typical $\omega$ when $d=1$ and $r>4$ \cite{SS}, pure point spectrum for typical $\omega$ when $d=1$ and $r>8$ (conjectured to be $r>2$) under some conditions on the density of distribution \cite{JM1}. 
We hope Theorem \ref{samain} could serve as a step towards the proof of the conjecture, as well as establishing localization in more general setting.

\

The second part of this note is devoted to extending some basic results which were previously known for random Schr\"odinger operators to long range self-adjoint cases.

Let $\Gamma$ be a subset of $\Z^d$ such that $\Gamma$ and $-\Gamma$ form a partition of $\Z^d\setminus \{0\}$ in the sense that $\Gamma \bigcap (-\Gamma)=\emptyset$ and $\Gamma\bigcup (-\Gamma)=\Z^d\setminus \{0\}$. 
In our definition of the operator, we will only specify the potentials $a_n$'s on $\Gamma$, the potentials on $-\Gamma$ will be defined automatically due to self-adjointness of the operator.
Denote $\Gamma_0=\{0\}\bigcup \Gamma$.
Let $\{\gamma^{(i)}\}_{i\in \Gamma}$ be a sequence of compactly supported probability measures on $\mathbb{C}$ and 
$\gamma^{(0)}$ be a probability measure on $\R$. 
Similar to Theorem \ref{samain} we do not make compact support assumption on $\gamma^{(0)}$. 
This enables us to cover unbounded operators.
Let $\mathrm{d}\kappa=\bigtimes_{i\in {\Gamma}_0} (\mathrm{d}\gamma^{(i)})^{\Z^d}$ and 
$\Omega=\bigtimes_{i\in \Gamma_0} (\supp\gamma^{(i)})^{\Z^d}
=\{\omega=(\omega^{(i)})_{i\in \Gamma_0}\ |\ \omega^{(i)}=(\omega^{(i)}_j)_{j\in \Z^d},\ \omega^{(i)}_j \in \mathrm{supp}(\gamma^{(i)})\}$. 
For any $n\in \Z^d$, we define the translations $\tilde{T}^n$ on $\Omega$ as 
$\tilde{T}^n \omega:=(T^n\omega^{(i)})_{i\in \Gamma_0}$, where $(T^n\omega^{(i)})_j=\omega^{(i)}_{j+n}$.

Consider long range self-adjoint operator $H_{\omega}$ on $l^2(\Z^d)$:
\begin{align}\label{defHomega}
(H_{\omega} u)(n)=\sum_{j\in \Z^d}a_{j}(\tilde{T}^n\omega )u(n-j),
\end{align}
where $a_k(\omega)=\omega^{(k)}_0$ for $k\in \Gamma_0$, and $a_{k}(\omega)=\overline{\omega^{(-k)}_{-k}}$ for $k\in (-\Gamma_0)$.
In particular, for $k=0$, $a_k(\omega)\in \R$.
We will further assume $\gamma^{(k)}$ have shrinking supports, namely, 
\begin{align}\label{supp}
\supp\gamma^{(k)}\subseteq B(C(1+\|k\|)^{-r})\ \ \mathrm{for}\ \ k\in \Gamma,
\end{align} 
for some constant $C>0$ and $r>\frac{d}{2}$, with $B(M)$ being the ball in $\mathbb{C}$ centered at $0$ with radius $M$. 

Following a standard argument, one can show that there exists a non-random $\Sigma$, 
such that $\sigma(H_{\omega})=\Sigma$ for $\kappa-$almost every $\omega\in\Omega$ \cite{KM2}.

We say a set $\Omega_1$ is a dense subset of $\Omega$, if for any $\omega \in \Omega$, $0<\xi\in \R$ and any finite sets $\bigwedge_1\subset \Gamma_0$ and $\bigwedge_2\subset \Z^d$, there exists $\tilde{\omega}\in \Omega_1$ such that $|\tilde{\omega}^{(i)}_j-\omega^{(i)}_j|<\xi$ for any $i\in \bigwedge_1$ and $j\in \bigwedge_2$.
It is clear that any subset $\Omega_0\subseteq \Omega$ with $\kappa(\Omega_0)=1$ is dense in $\Omega$.

The following result is similar to Theorem 3 \cite{KM}.
\begin{thm}\label{densespec}
Let $\Omega_1$ be a dense subset of $\Omega$, then
\begin{align*}
\Sigma=\overline{\bigcup_{\omega\in\Omega_1}\sigma(H_{\omega})}
\end{align*}
\end{thm}
\

\

This theorem has a direct corollary which implies the almost sure spectrum is determined by the spectra of periodic operators. 
We say $\omega=(\omega^{(i)})_{i\in \Gamma_0}\in \Omega$ is $p-$periodic if $\omega^{(i)}_j=\omega^{(i)}_{j+p}$ for any $i\in \Gamma_0$ and $j\in \Z^d$.
\begin{cor}\label{pper}
\begin{align*}
\Sigma=\overline{\bigcup_{p\in \Z_+}\ \ \bigcup_{p-\mathrm{periodic}\, \omega}\sigma(H_{\omega})}
\end{align*}
\end{cor}
This can be viewed as an extension of Theorem 3.9 \cite{Invitation} for Schr\"odinger operators. 

Contrast to the general long range case, if we focus on one-dimensional random Jacobi matrices, it turns out we could get a better result than Corollary \ref{pper}. 
Let $\gamma$ be a compactly supported measure on $\R$.
Let $\Gamma=\{1\}$ and $\Gamma_0=\{0,1\}$. 
Let $d\tilde{\kappa}=\bigtimes_{i\in \Gamma_0}(\mathrm{d}\gamma)^{\Z}$ and $\tilde{\Omega}=\bigtimes_{i\in \Gamma_0}(\supp{\gamma})^{\Z}$.
For any $\omega\in \tilde{\Omega}$, we consider
\begin{align}\label{defJacobiH}
(H_\omega u)(n)=a(T^n\omega^{(1)})u(n+1)+a(T^{n-1}\omega^{(1)})u(n-1)+b(T^n\omega^{(0)})u(n)
\end{align}
where $a(\omega^{(1)})=\omega_0^{(1)}$, $b(\omega^{(0)})=\omega_0^{(0)}$ and $(T\omega^{(i)})_n=\omega^{(i)}_{n+1}$ for $i=0,1$. We have
\begin{thm}\label{2per} Let $H_{\omega}$ be a random Jacobi matrix defined as in (\ref{defJacobiH}). We have
\begin{align*}
\Sigma=\bigcup_{1-\mathrm{periodic}\, \omega}\sigma(H_{\omega}).
\end{align*}
\end{thm}

We organize the note as follows: a key lemma and proof of part (a) of Theorem \ref{main} will be presented in Section 2, proofs of Theorems \ref{densespec} and \ref{2per} will be given in Section 3.
The proofs of parts (b), (c) of Theorem \ref{main} will be included in the appendix.

\section{Key lemma and the proof of Theorem \ref{main}}
In order to prove Theorem \ref{main}, we need the following lemma.
\begin{lemma}\label{key}
Let $H$ be defined as in (\ref{defnormalH}) with $r>(2d+\epsilon)\frac{q}{q-1}$ for some $\epsilon> 0$ and $2\leq q\in \Z^{+}$.
Let $z\in G_{\epsilon}$ and $\phi$ be a corresponding generalized eigenfunction. Let 
\begin{align}\label{defphiN}
 \phi_N(j)=\left\{
\begin{array}{lcl}
\phi(j)     &      & \|j\|\leq N,\\
0            &      & otherwise.
\end{array} \right.
\end{align}
and $\Phi_N=(H-z)\phi_N$. We have
\begin{align}\label{eq:Lemma21}
\liminf_{L\rightarrow \infty}\frac{\|\Phi_{L^q}\|_{l^2}}{\|\phi_{L^q}\|_{l^2}}=0.
\end{align}
\end{lemma}

We will first give the proof of Theorem \ref{main} based on Lemma \ref{key}.
\subsection*{Proof of Theorem \ref{main}}
We only prove part (a) here, the proofs of part (b) and (c) are standard, we include them in the appendix for reader's convenience.

Since $r>2d+\epsilon$, we can choose $q\in \Z^+$ large enough such that $r>(2d+\epsilon)\frac{q}{q-1}$.
Then Lemma \ref{key} directly implies $\dist{(z, \sigma(H))}=0$. $\hfill{} \Box$

\subsection*{Proof of Lemma \ref{key}}
For simplicity, we denote $a_0^j-z$ by $a_0^j$ and let $S_r:=(\sum_{m\in \Z^d}(1+\|m\|)^{-r})^2$.  
We have $\Phi_N(n)=\sum\limits_{\|j\|\leq N}a_{n-j}^j \phi(j)$.
In view of the numerator of \eqref{eq:Lemma21}, we have
\begin{align}\label{1}
\sum_{n\in \Z^d} |\Phi_{L^q}(n)|^2
=&\sum_{n\in \Z^d} |\sum_{\|j\|\leq L^q}a_{n-j}^j \phi(j)|^2 \notag\\
=&\sum_{\|n\|> L^q}|\sum_{\|j\|\leq L^q, \|j-n\|> L^{q-1}}a_{n-j}^j\phi(j)+\sum_{\|j\|\leq L^q, \|j-n\|\leq L^{q-1}}a_{n-j}^j\phi(j)|^2 \notag\\
 &+\sum_{\|n\|\leq L^q}|\sum_{\|j\|>L^q, \|j-n\|>L^{q-1}}a_{n-j}^j\phi(j)+\sum_{\|j\|>L^q, \|j-n\|\leq L^{q-1}}a_{n-j}^j\phi(j)|^2 \notag\\
\leq &2\sum_{\|n\|> L^q}\left(|\sum_{\|j\|\leq L^q, \|j-n\|> L^{q-1}}a_{n-j}^j\phi(j)|^2+|\sum_{\|j\|\leq L^q, \|j-n\|\leq L^{q-1}}a_{n-j}^j\phi(j)|^2\right) \notag\\
       &+2\sum_{\|n\|\leq L^q}\left(|\sum_{\|j\|>L^q, \|j-n\|>L^{q-1}}a_{n-j}^j\phi(j)|^2+|\sum_{\|j\|> L^q, \|j-n\|\leq L^{q-1}}a_{n-j}^j\phi(j)|^2 \right) \notag\\
=:&2(\Sigma_1+\Sigma_2+\Sigma_3+\Sigma_4).
\end{align}
Next, we are going to estimate the four terms separately.

\subsection*{Estimate of $\Sigma_1$}
By Cauchy-Schwarz inequality and our assumption on $a_{n-j}$, we have
\begin{align*}
       \Sigma_1=\sum_{\|n\|>L^q}|\sum_{\|j\|<L^q, \|j-n\|>L^{q-1}}a_{n-j}^j\phi(j)|^2 
\leq &\|\phi_{L^q}\|_{l^2}^2 \sum_{\|n\|>L^q}\ \ \sum_{\|j\|<L^q, \|j-n\|>L^{q-1}}|a_{n-j}^j|^2 \notag\\
\leq &\|\phi_{L^q}\|_{l^2}^2 \sum_{\|n\|>L^q}\ \ \sum_{\|j\|<L^q, \|j-n\|>L^{q-1}}(1+\|n-j\|)^{-2r}.
\end{align*}
Dividing the sum in $n$ into two parts and estimate each part separately, we have
\begin{align}\label{Sigma1app}
\Sigma_1 \leq &\|\phi_{L^q}\|_{l^2}^2 \left( \sum_{L^q+L^{q-1}\geq \|n\|>L^q}\ \ \sum_{\|j\|<L^q, \|n-j\|>L^{q-1}}(1+\|n-j\|)^{-2r} \right. \notag\\
&\left. \qquad\qquad +\sum_{\|n\|>L^q+L^{q-1}}\ \ \sum_{\|j\|<L^q}(1+\|n-j\|)^{-2r} \right) \notag\\
\leq & C\|\phi_{L^q}\|_{l^2}^2 \left(L^{q(2d-2r)+2r-1}+L^{qd}\int_{L^q+L^{q-1}}^{\infty} \frac{x^{d-1}}{(x-L^q)^{2r}}\ \mathrm{d}x\right) \notag\\
\leq & C\|\phi_{L^q}\|_{l^2}^2 L^{q(2d-2r)+2r-1}.
\end{align}

\subsection*{Estimate of $\Sigma_2$}
\begin{align}
  \Sigma_2=&\sum_{\|n\|> L^q}|\sum_{\|j\|\leq L^q, \|j-n\|\leq L^{q-1}}a_{n-j}^j\phi(j)|^2\notag\\
\leq &\sum_{L^q+L^{q-1}\geq \|n\|>L^q}\left( \sum_{\|m\|\leq L^{q-1}}\chi_{\|n-m\|\leq L^q}\ |a_{m}^{n-m}\phi(n-m)| \right)^2\notag\\
\leq &\left\lbrace \sum_{\|m\|\leq L^{q-1}} (1+\|m\|)^{-r}  \left(\sum_{L^q+L^{q-1}\geq \|n\|>L^q} \chi_{\|n-m\|\leq L^q} |\phi(n-m)|^2\right)^{\frac{1}{2}} \right\rbrace^2\label{Sigma2'}\\
\leq &S_r\ (\|\phi_{(L+1)^q}\|_{l^2}^2-\|\phi_{(L-1)^q}\|_{l^2}^2).\label{Sigma2}
\end{align}
We used Minkowski's inequality in \eqref{Sigma2'}.

\subsection*{Estimate of $\Sigma_3$}
By Cauchy-Schwarz inequality and Minkowski's inequality, we have
\begin{align*}
       \Sigma_3=&\sum_{\|n\|\leq L^q}|\sum_{\|j\|>L^q, \|j-n\|>L^{q-1}}a_{n-j}^j\phi(j)|^2\\
\leq &\sum_{\|n\|\leq L^q}(\sum_{\|j\|>L^q, \|j-n\|>L^{q-1}}|a_{n-j}^j\phi(j)|)^2\\
\leq &\left\lbrace \sum_{\|j\|>L^q} |\phi(j)| \left( \sum_{\|n\|\leq L^q} \chi_{\|j-n\|>L^{q-1}}\ |a_{n-j}^j|^2 \right)^{\frac{1}{2}} \right\rbrace^2.
\end{align*}
Using our assumption on $a_{n-j}$, and the fact that $\phi$ is an $\epsilon$-generalized eigenfunction, we have
\begin{align}
\Sigma_3\leq &\left\lbrace \sum_{L^q+L^{q-1}\geq \|j\|>L^q} \|j\|^{\frac{d}{2}+\epsilon} \left( \sum_{\|n\|\leq L^q} \chi_{\|j-n\|>L^{q-1}}\ (1+\|n-j\|)^{-2r} \right)^{\frac{1}{2}}\right. \notag\\
       &\left.+\sum_{\|j\|>L^q+L^{q-1}} \|j\|^{\frac{d}{2}+\epsilon} \left(\sum_{\|n\|\leq L^q} (1+\|n-j\|)^{-2r}\right)^{\frac{1}{2}} \right\rbrace^2 \notag\\
\leq &\left\lbrace CL^{q(2d-r)+r+\epsilon q}+\sum_{\|j\|>L^q+L^{q-1}} \|j\|^{\frac{d}{2}+\epsilon} \left(\sum_{\|n\|\leq L^q} (1+\|n-j\|)^{-2r}\right)^{\frac{1}{2}} \right\rbrace^2 \notag\\
\leq &\left\lbrace CL^{q(2d-r)+r+\epsilon q}+L^{\frac{dq}{2}} \sum_{\|j\|>L^q+L^{q-1}}(1+\|j\|)^{\frac{d}{2}+\epsilon}\ (\|j\|-L^q)^{-r} \right\rbrace^2 \notag\\
\leq &\left\lbrace CL^{q(2d-r)+r+\epsilon q}+C L^{\frac{dq}{2}}\int_{L^q+L^{q-1}}^{\infty}\frac{x^{\frac{3d}{2}+\epsilon-1}}{(x-L^q)^r}\ \mathrm{d}x\right\rbrace^2 \notag\\
\leq &(CL^{q(2d-r)+r+\epsilon q})^{\frac{1}{2}}.   \label{Sigma3app}
\end{align}

\subsection*{Estimate of $\Sigma_4$}
\begin{align}
  \Sigma_4=&\sum_{\|n\|\leq L^q}|\sum_{\|j\|> L^q, \|j-n\|\leq L^{q-1}}a_{n-j}^j\phi(j)|^2\notag\\
\leq &\sum_{\|n\|\leq L^q} \left(\sum_{\|m\|\leq L^{q-1}}\chi_{\|n-m\|> L^q}|a_{m}^{n-m}\phi(n-m)|\right)^2\notag\\
\leq &\left\lbrace \sum_{\|m\|\leq L^{q-1}} (1+\|m\|)^{-r} \left(\sum_{\|n\|\leq L^q} \chi_{\|n-m\|> L^q} |\phi(n-m)|^2\right)^{\frac{1}{2}}\right\rbrace^2\label{Sigma4'}\\
\leq &S_r\ (\|\phi_{(L+1)^q}\|_{l^2}^2-\|\phi_{(L-1)^q}\|_{l^2}^2).\label{Sigma4}
\end{align}
We used Minkowski's inequality in \eqref{Sigma4'}.

Combining the estimates (\ref{Sigma1app}), (\ref{Sigma2}), (\ref{Sigma3app}), (\ref{Sigma4}), we obtain
\begin{align}
&\liminf_{L\rightarrow\infty}\frac{\|\Phi_{L^q}\|_{l^2}^2}{\|\phi_{L^q}\|_{l^2}^2}\notag\\
\leq
&\liminf_{L\rightarrow\infty}\left(2 S_r\ \frac{\|\phi_{(L+1)^q}\|_{l^2}^2-\|\phi_{(L-1)^q}\|_{l^2}^2}{\|\phi_{L^q}\|_{l^2}^2}+
C L^{q(2d-2r)+2r-1}+ \frac{(CL^{q(2d-r)+r+\epsilon q})^{\frac{1}{2}}}{\|\phi_{L^q}\|_{l^2}^2}\right).\label{eq:combine}
\end{align}
Note that by our choice of $r>(2d+\epsilon)q/(q-1)$, we have 
$$\lim_{L\rightarrow\infty}L^{q(2d-2r)+2r-1}=\lim_{L\rightarrow\infty}L^{q(2d-r)+r+\epsilon q}=0,$$
and
$$S_r=\sum_{m\in \Z^d}(1+\|m\|)^{-r}<\infty.$$
Also note that by our definition of $\epsilon$-generalized eigenfunction, we have $\phi(0)=1$, hence $\|\phi_{L^q}\|_{l^2}\geq 1$.
\linebreak Therefore \eqref{eq:combine} yields
\begin{align}
\liminf_{L\rightarrow\infty}\frac{\|\Phi_{L^q}\|_{l^2}^2}{\|\phi_{L^q}\|_{l^2}^2}\leq C\liminf_{L\rightarrow\infty}\frac{\|\phi_{(L+1)^q}\|_{l^2}^2-\|\phi_{(L-1)^q}\|_{l^2}^2}{\|\phi_{L^q}\|_{l^2}^2}.
\end{align}
Then  $\liminf_{L\rightarrow\infty}\frac{\|\Phi_{L^q}\|_{l^2}^2}{\|\phi_{L^q}\|_{l^2}^2}>\kappa>0$ would imply that
$\|\phi_{(L+1)^q}\|_{l^2}^2\geq (1+\frac{\kappa}{C})\|\phi_{(L-1)^q}\|_{l^2}^2$, which leads to exponential growth of $\phi$, contradiction.
$\hfill{} \Box$

\section{Proofs of Theorems \ref{densespec} and \ref{2per}}
\subsection{Proof of Theorem \ref{densespec}}
Let $\Omega_0$ be a shift invariant set with $\kappa(\Omega_0)=1$ such that $\sigma(H_{\omega})=\Sigma$ for every $\omega\in\Omega_0$. 
Then given any $\omega\in\Omega$,  $\xi>0$, and any finite sets $\bigwedge_1 \subset \Gamma_0$, $\bigwedge_2 \subset \Z^d$, 
there exists $\tilde{\omega}\in\Omega_0$ such that $|\omega^{(j)}_n-\tilde{\omega}^{(j)}_n|<\xi$ for any $j\in\bigwedge_1$ and $n\in\bigwedge_2$.

\subsection*{The $``\supseteq"$ direction}
Take $\omega\in \Omega_1$ and $E\in \sigma(H_{\omega})$. By Weyl's criterion, for any $L>0$ there exists $\phi^{(L)}\in l^2(\Z^d)$ such that 
$\|(H_{\omega}-E)\phi^{(L)}\|<\frac{1}{L}\|\phi^{(L)}\|$.

First, we show such $\phi^{(L)}$ can be taken with compact support. 
This is standard if $H_{\omega}$ is a bounded operator, but since $H_{\omega}$ could be unbounded here, we will work out the details.

Let us consider a cut-off function $\phi^{(L)}_k(n)=\chi_{\|j\|\leq k}(n) \phi^{(L)}(n)$. 
We split $H_{\omega}$ into two parts $H_{\omega}^{0}+H_{\omega}^{1}$, where $H_{\omega}^0$ contains only the diagonal multiplication, namely $(H_{\omega}^{0}u)(n)=a_0(\tilde{T}^n \omega)u(n)$, and $H_{\omega}^{1}$ is the off-diagonal part.
Clearly, $H_{\omega}^{1}$ is a bounded operator while $H_{\omega}^{0}$ could be unbounded.
Now let us consider 
\begin{align*}
  (H_{\omega}-E)(\phi^{(L)}-\phi^{(L)}_k)=(H_{\omega}^{1}-E)(\phi^{(L)}-\phi^{(L)}_k)+H_{\omega}^{0}(\phi^{(L)}-\phi^{(L)}_k).
\end{align*}
We know $\|(H_{\omega}^{1}-E)(\phi^{(L)}-\phi^{(L)}_k)\|_{l^2}\leq \|(H_{\omega}^{1}-E)\|\ \|(\phi^{(L)}-\phi^{(L)}_k)\|_{l^2}\rightarrow 0$ as $k\rightarrow\infty$.
Also since $\|H_{\omega}^{0}\phi^{(L)}\|\leq \|(H_{\omega}^{1}-E)\phi^{(L)}\|+\frac{1}{L}\|\phi^{(L)}\|<\infty$, 
we have $\|H_{\omega}^{0}(\phi^{(L)}-\phi^{(L)}_k)\|_{l^2}\rightarrow 0$, as $k\rightarrow\infty$. 
Thus $\|(H_{\omega}-E)(\phi^{(L)}-\phi^{(L)}_k)\|\rightarrow 0$ as $k\rightarrow \infty$. 
Taking $k_L$ large enough, we may assume 
\begin{align}
\|(H_{\omega}-E)\phi^{(L)}_{k_L}\|<\frac{2}{L}\|\phi^{(L)}_{k_L}\|.
\end{align}
For simplicity, we still denote $\phi^{(L)}_{k_L}$ by $\phi^{(L)}$ and we assume it is supported on a compact set $\{\|n\|\leq K\}$.

\

Now, we take an integer $m$ large enough such that $(m-1)^{2r-d}>L^2 K^{2d-2r}$, this choice is guaranteed by $r>d/2$. 
Let $\bigwedge_1=\Gamma_0\bigcap \{\|j\|\leq (m+1)K\}$ and $\bigwedge_2=\{\|n\|\leq mK\}$.
There exists $\tilde{\omega}\in \Omega_0$ such that $|\omega^{(j)}_n-\tilde{\omega}^{(j)}_n|<\frac{1}{LK^d m^{d/2}}$ for any $j\in \bigwedge_1$ and $n\in\bigwedge_2$.
Now let us consider
\begin{align*}
  &\|(H_{\tilde{\omega}}-H_{{\omega}})\phi^{(L)}\|_{l^2}^2\\
=&\sum_{n\in \Z^d} |\sum_{\|j\|\leq K} (a_{n-j} (\tilde{T}^n \tilde{\omega})-a_{n-j} (\tilde{T}^n\omega)) \phi^{(L)}(j)|^2\\
\leq &\|\phi^{(L)}\|_{l^2}^2 \sum_{n\in \Z^d}\ \sum_{\|j\|\leq K} |a_{n-j} (\tilde{T}^n \tilde{\omega})-a_{n-j} (\tilde{T}^n \omega)|^2\\
=& \|\phi^{(L)}\|_{l^2}^2 \sum_{\|j\|\leq K}\ (\sum_{\|n\|>mK}+\sum_{\|n\|\leq mK, n-j\in \Gamma_0}+\sum_{\|n\|\leq mK, n-j\in -\Gamma})
|a_{n-j} (\tilde{T}^n \tilde{\omega})-a_{n-j} (\tilde{T}^n \omega)|^2\\
\leq &\|\phi^{(L)}\|_{l^2}^2 \sum_{\|j\|\leq K}\ (\sum_{\|n\|>mK} 2(1+\|n-j\|)^{-2r}+\sum_{\|n\|\leq mK, n-j\in \Gamma_0}|\tilde{\omega}_n^{(n-j)}-\omega_n^{(n-j)}|^2\\
&\ \ \ \ \ \ \ \ \ \ \ \ \ \ \ \ \ \ \ \ \ \ \ \ \ \ \ \ \ \ \ \ \ \ \ \ \ \ \ \ \ \ \ \ \ \ \ \ \ \ \ \ \ \ \ \ \ \ \ \ \ \ \ \ \ \ \ \ \ \ \ \ \ +\sum_{\|n\|\leq mK, n-j\in -\Gamma}|\tilde{\omega}^{(j-n)}_n-\omega^{(j-n)}_n|^2)\\
\leq &C\|\phi^{(L)}\|_{l^2}^2 (\sum_{\|n\|\geq mK} \sum_{\|j\|\leq K} \|n-j\|^{-2r}+K^{2d}m^d\frac{1}{L^2 K^{2d} m^d})\\
\leq &C\|\phi^{(L)}\|_{l^2}^2 ((m-1)^{d-2r}K^{2d-2r}+\frac{1}{L^2})\\
\leq &\frac{C}{L^2}\|\phi^{(L)}\|_{l^2}^2
\end{align*}
Thus $\|(H_{\tilde{\omega}}-E)\phi^{(L)}\|\leq \frac{\sqrt{C}+1}{L}\|\phi^{(L)}\|$. Taking $L\rightarrow\infty$, we get $E\in \sigma(H_{\tilde{\omega}})=\Sigma$.

\subsection*{The $``\subseteq"$ direction}
Note that since $\Omega_1$ is dense in $\Omega$, the proof follows from that of $``\supseteq"$ by interchanging the roles of $\Omega_0$ and $\Omega_1$.$\hfill{} \Box$

\subsection{Proof of Theorem \ref{2per}}
Let $H_\omega=H_\omega^1+H_\omega^2$, where $(H_\omega^1 u)(n)=a(T^n\omega^{(1)})u(n+1)+a(T^{n-1}\omega^{(1)})u(n-1)$ 
and $(H_\omega^2 u)(n)=b(T^n\omega^{(0)})u(n)$. Let $\Sigma_1=\sigma(H_\omega^1)$ a.s.,  $\Sigma_2=\sigma(H_\omega^2)$ a.s. and $M=\sup\{|\alpha|, \alpha\in\supp\gamma\}$.
Clearly, 
\begin{align}
\bigcup_{1-\mathrm{periodic}\, \omega}\sigma(H_{\omega})=[-2M, 2M]+\supp{\gamma}.
\end{align}

\subsection*{The $``\subseteq"$ direction}
$\|H_\omega^1\|\leq 2M$, we have $\Sigma_1\subset [-2M, 2M]$. 
This immediately implies $\Sigma\subseteq\Sigma_1+\Sigma_2\subseteq [-2M, 2M]+{\supp}\gamma$. 
\subsection*{The $``\supseteq"$ direction}
Let $\Omega_0$ be the full $\tilde{\kappa}$-measure set so that $\sigma(H_{\omega})=\Sigma$ for any $\omega\in \Omega_0$.
Since $M=\sup\{|\alpha|, \alpha\in\supp\gamma\}$, for any $\beta\in (-2M, 2M)$, there exists a set $\mathcal{F}$ with $\gamma(\mathcal{F})>0$ such that $\frac{1}{2}|\beta|<\inf\{|\alpha|, \alpha\in \mathcal{F}\}$. 
Then taking any $\xi\in \supp\gamma$, there exists a sequence of 1-periodic $\omega(n)=(\omega^{(1)}(n), \omega^{(0)}(n))\in \Omega_0$ such that 
for any $k\in \Z$, $\omega_k^{(1)}(n)\equiv \alpha_n$ and $\omega_k^{(2)}(n)\equiv \xi_n$, with $\frac{1}{2}|\beta|<\alpha_n$ and $\xi_n\rightarrow \xi$ as $n\rightarrow\infty$. 
Clearly, this implies $\beta+\xi\in \bigcup_{n} \sigma(H_{\omega(n)})=\Sigma$. Thus we have $(-2M, 2M)+\supp\gamma\subseteq \Sigma$, which implies the desired result after taking closure. $\hfill{} \Box$

\section*{\\Appendix}

\subsection*{Part (b) or Theorem \ref{main}}
The proof is standard, we present it here for readers' convenience and completeness.

Since $\mu_{n,m}\ll \mu$, there exists a $F_{n,E}(m)=\frac{d\mu_{n,m}}{d\mu}(E)$ defined for $\mu$-a.e. $E$.
We will show that for any fixed $n$, for $\mu$-a.e. $E$, $F_{n,E}(m)$ is an $\epsilon$-generalized eigenfunction, namely
\begin{enumerate}
\item $|F_{n,E}(m)|\leq C(n) (1+\|m\|)^{\frac{d}{2}+\epsilon}$,
\item $F_{n,E}(m)$ is a solution to $Hu=Eu$.
\end{enumerate}
\begin{proof}

(1). For any Borel set $B\subset\mathbb{C}$, 
$$
|\frac{1}{\mu(B)}\int_{B}F_{n,E}(m)\ d\mu(E)|
=|\frac{\mu_{n,m}(B)}{\mu(B)}|
\leq \frac{\sqrt{\mu_n(B)\mu_m(B)}}{\mu(B)}|
\leq \frac{1}{\sqrt{\lambda_n\lambda_m}}
\leq C (1+\|n\|)^{\frac{d}{2}+\epsilon}(1+\|m\|)^{\frac{d}{2}+\epsilon}.
$$
Lebesgue differentiation theorem implies that 
$|F_{n,E}(m)|\leq C(n) (1+\|m\|)^{\frac{d}{2}+\epsilon}$ for $\mu$-a.e. $E$.

(2). $((H-E)F_{n,E})(m)=0$ for $\mu$-a.e. $E$, is equivalent to 
\begin{align}\label{b2}
\int_{\sigma(H)} (HF_{n,E})(m)\, g(E)\ d\mu(E)=\int_{\sigma(H)} F_{n,E}(m)\, E\, g(E)\ d\mu(E)
\end{align}
for any compactly supported continuous function $g$ on $\sigma(H)$. 

First, note that $HF_{n,E}(m)$ is well defined for $\mu$-a.e. $E$. 
Indeed, by definition,
\begin{align}\label{series}
HF_{n,E}(m)=\sum_{k\in \Z^d} a^m_{m-k} F_{n,E}(k).
\end{align}
On one hand, for any fixed $n$, by our estimates in (1), $|F_{n,E}(k)|$ is bounded above by $C(1+\|k\|)^{\frac{d}{2}+\epsilon}$.
On the other hand, $a^m_{m-k}$ decays like $C(1+\|m-k\|)^{-r}$ with $r>2d$.
Hence for any fixed $m$, the convergence of the series in \eqref{series} is verified.

To prove \eqref{b2}, note that for any compactly supported continuous function $g$, $g(H)$ is a bounded operator.
We have
\begin{align*}
\int_{\sigma(H)} F_{n,E}(m)\, E\, g(E)\ d\mu(E)
=&(e_n, Hg(H)e_m)\\
=&(g(H)^{*}e_n, He_m)\\
=&(g(H)^{*}e_n, \sum_{k\in \Z^d}a_{k-m}^me_k)\\
=&(e_n, \sum_{k\in \Z^d}a_{k-m}^m g(H)e_k)\\
=&\int_{\sigma(H)}\sum_{k\in \Z^d}a_{k-m}^m F_{n,E}(k)g(E)\ d\mu(E)\\
=&\int_{\sigma(H)} (HF_{n,E})(m)\, g(E)\ d\mu(E).
\end{align*}
This proves \ref{b2}.

\end{proof}
$\hfill{} \Box$

\subsection*{Part (c) of Theorem 1.1}
Part (c) follows directly from (a), (b) and the fact that spectrum of $H$ is the smallest closed set which supports every spectral measure of $H$. $\hfill{} \Box$

\section*{Acknowledgement}
I would like to thank Svetlana Jitomirskaya for suggesting this problem and for many useful discussions.
This research was partially supported by the NSF DMS-1401204.

\bibliographystyle{amsplain}

\end{document}